\documentclass[a4paper,12pt]{article}
\usepackage{amssymb}

\usepackage{amsmath}
\usepackage{a4}

\input{tcilatex}

\begin{document}

\title{A homotopy aspect of representation theory}
\author{Haibao Duan \\
Institute of Mathematics, Chinese Academy of Sciences\\
Beijing 100080, dhb@math.ac}
\date{ \ \ \ \ \ \ }
\maketitle

\begin{abstract}
We bring a linkage from representation theory of Lie groups to homotpy
theory for maps between flag manifolds. As applications we derive from
representation theory abundant families of distinguishable homotopy classes
of maps between flag manifolds.

\begin{description}
\item \textsl{2000 Mathematical Subject Classification: }57R22 (22E46).
\end{description}
\end{abstract}

\begin{center}
\textbf{1. Introduction }
\end{center}

A current topic in algebraic geometry is the relationship between
representation theory of Lie groups and Schubert calculus in flag manifolds
(cf. [LG, Chapter 4], [K, VIII], [T]). It seems not aware in the existing
literature that there is a direct connection between representation theory
and homotopy theory for maps between flag manifolds, implicitly contained in
the classical works of Borel and Hirzebruch [BH]. This connection may be
summarized as ``\textsl{character factors through homotopy}'' and is
presented in the Theorem in \S 5. Applying this idea we obtain abundant
families of maps between flag manifolds whose homotopy classes are
distinguishable by the characters of representations (Corollary 1--3 in \S
5).

We begin with a problem concerning homotopy classification for maps between
flag manifolds (\S 2), followed by a brief introduction to what we need from
representation theory (\S 3,\S 4). The main results are stated and
established respectively in \S 5 and \S 6. The final section \S 7 is devoted
to further comments.

\begin{center}
\textbf{2. A realization problem in homotopy theory }
\end{center}

For a compact connected semi--simple Lie group $G$ with a maximal torus $%
T\subset G$, write $G/T$ for the flag manifold of $G$ modulo $T$. If $G$ is
the unitary group $U(n)$ of order $n$ and if $T^{n}\subset U(n)$ is the
diagonal subgroup, we write $F(n)$ instead of $U(n)/T^{n}$.

The primary objective of our study is the set $[G/T,F(n)]$ of all homotopy
classes of maps from $G/T$ to $F(n)$ (cf. \S 7.1 for references on the
topic). For this purpose one may examine the correspondence that sends a map
$f$ to the induced cohomology homomorphism $f^{\ast }$.

\begin{enumerate}
\item[(A)] \qquad \qquad \noindent $r:[G/T,F(n)]\rightarrow
Hom[H^{2}(F(n)),H^{2}(G/T)]$

\noindent \qquad (resp. $r_{1}:[G/T,F(n)]\rightarrow Hom[H^{\ast
}(F(n)),H^{\ast }(G/T)]$),
\end{enumerate}

\noindent where $Hom[H^{\ast }(F(n)),H^{\ast }(G/T)]$ is the set of all ring
maps from $H^{\ast }(F(n))$ to $H^{\ast }(G/T)$. It follows from rational
homotopy theory [GH] that

\begin{quote}
\textsl{the }$r_{1}$\textsl{\ is finite to one;}
\end{quote}

\noindent We recall also from Bott and Samelson [BS] that

\begin{quote}
\textsl{the ring }$H^{\ast }(G/T)$\textsl{\ is torsion free and over
rationals generated multiplicatively by elements from }$H^{2}(G/T)$\textsl{.}
\end{quote}

\noindent For these reasons a solution to the following problem is an
important step towards understanding the set $[G/T,F(n)]$.

\begin{quote}
\textbf{Problem.} \textsl{Which homomorphisms }$h:H^{2}(F(n))\rightarrow
H^{2}(G/T)$\textsl{\ are induced by maps }$G/T\rightarrow F(n)$\textsl{?}
\end{quote}

\begin{center}
\textbf{3. Notions}
\end{center}

For a compact connected Lie group $G$ with a fixed maximal torus $T$,
consider the transgression

\begin{enumerate}
\item[(B)] $\qquad \qquad \tau _{G}:H^{1}(T)=Hom(T,S^{1})\rightarrow
H^{2}(G/T)$
\end{enumerate}

\noindent in the fibration $T\hookrightarrow G\rightarrow G/T$ [BH, \S
10.1], where $S^{1}$ is the circle group, and where as indicated the $1$%
--dimensional cohomology $H^{1}(T)$ is canonically identified with $%
Hom(T,S^{1})$, the set of $1$--dimensional unitary representation of $T$
(cf. Lemma 6 in \S 6).

If $G=U(n)$ is the unitary group, the diagonal subgroup $T^{n}\subset U(n)$
has the factorization $T^{n}=diag\{e^{i\theta _{1}},\cdots ,e^{i\theta
_{n}}\}$. Let $t_{k}^{\prime }\in H^{1}(T^{n})=Hom(T^{n},S^{1})$ be the
projection onto the $k^{th}$ factor, $1\leq k\leq n$, and put $t_{k}=\tau
_{U(n)}(t_{k}^{\prime })\in H^{2}(F(n))$.

\bigskip

\textbf{Lemma 1.} \textsl{In }$H^{2}(F(n))$ \textsl{one has }$t_{1}+\cdots
+t_{n}=0$\textsl{. Further,} \textsl{the classes }$t_{1},\cdots ,t_{n-1}$%
\textsl{\ form a basis of }$H^{2}(F(n))$\textsl{.}

\textbf{Proof.} This follows from the description of the ring $H^{\ast
}(F(n))$ due to Borel [B]: $H^{\ast }(F(n))=\mathbb{Z}[t_{1},\cdots
,t_{n}]/<e_{i}\mid 1\leq i\leq n>$, where the $e_{i}\in \mathbb{Z}%
[t_{1},\cdots ,t_{n}]$ is the $i^{th}$ elementary symmetric function in the $%
t_{1},\cdots ,t_{n}$.$\square $

\bigskip

Consider next a compact connected semi--simple $G$ with a fixed maximal
torus $T$. Let $p:\widetilde{G}\rightarrow G$ be the universal cover of $G$
and let $\widetilde{T}\subset \widetilde{G}$ be the maximal torus that
corresponds to $T$ under $p$. Then $p$ induces an isomorphism $\widetilde{G}/%
\widetilde{T}\cong G/T$ between flag varieties. It is known that ([BH,\S
10.1], [DZZ,Theorem 2])

\bigskip

\textbf{Lemma 2.} \textsl{A set }$\{\omega _{i}\in Hom(\widetilde{T}%
,S^{1})\mid 1\leq i\leq m=\dim \widetilde{T}\}$\textsl{\ of fundamental
dominant weights of }$\widetilde{G}$ \textsl{relative to }$\widetilde{T}$%
\textsl{\ constitutes a basis of }$H^{1}(\widetilde{T})=Hom(\widetilde{T}%
,S^{1})$ \textsl{(cf. Lemma 6 in \S 6), and the transgression }$\tau _{%
\widetilde{G}}:H^{1}(\widetilde{T})\rightarrow H^{2}(G/T)$ \textsl{is an
isomorphism.}

\bigskip

Consequently, if we make no difference in notation between elements in $%
H^{1}(\widetilde{T})$ and their images under $\tau _{\widetilde{G}}$ in $%
H^{2}(G/T)$ (as in [BH]), we have

\bigskip

\textbf{Lemma 3.} \textsl{The set }$\{\omega _{i}\in H^{2}(G/T)\mid 1\leq
i\leq m\}$\textsl{\ of fundamental dominant weights of }$\widetilde{G}$%
\textsl{\ is a basis of }$H^{2}(G/T)$\textsl{.}

\bigskip

\textbf{Remark 1.} The classes $\omega _{i}\in H^{2}(G/T)$, $1\leq i\leq m$,
are of particular interests in the algebraic intersection theory on $G/T$.
With respect to the classical Schubert cell decomposition of $G/T$ [BGG],
they are precisely the \textsl{special Schubert classes} on $G/T$ (cf. [DZZ,
\S 5.2]).

\bigskip\

In view of Lemma 1 and 3, every homomorphism $h:H^{2}(F(n))\rightarrow
H^{2}(G/T)$ has a canonical numerical characterization as

\begin{enumerate}
\item[(C)] \qquad \qquad \noindent $h(t_{k})=a_{k,1}\omega _{1}+\cdots
+a_{k,m}\omega _{m}$, $1\leq k\leq n-1$,
\end{enumerate}

\noindent where $a_{k,i}\in \mathbb{Z}$. As a result, if we let $\mathbb{Z}%
^{+}$ be the set of non--negative integers, and let $\mathbb{Z}^{+}[\omega
_{1},\cdots ,\omega _{m},\rho ]$ be the semiring of polynomials in the
variables $\omega _{1},\cdots ,\omega _{m},\rho $ with coefficients in $%
\mathbb{Z}^{+}$ that is subject to the relation

\begin{enumerate}
\item[(D)] \qquad \qquad \noindent $\omega _{1}\cdots \omega _{m}\rho =1$,
\end{enumerate}

\noindent we may introduce a map

\begin{center}
$s:Hom(H^{2}(F(n)),H^{2}(G/T))\rightarrow \mathbb{Z}^{+}[\omega _{1},\cdots
,\omega _{m},\rho ]$ \
\end{center}

\noindent by the simple algorithm: if $h:H^{2}(F(n))\rightarrow H^{2}(G/T)$
is given as that in (C), we put

\begin{center}
$s(h)=\sum\limits_{1\leq k\leq n-1}$\noindent $\omega _{1}^{a_{k,1}}\cdots
\omega _{m}^{a_{k,m}}+$\noindent $\omega _{1}^{b_{1}}\cdots \omega
_{m}^{b_{m}}$,
\end{center}

\noindent where $b_{i}=-(a_{1,i}+\cdots +a_{n-1,i})$. We note that although
the $s(h)$ may appear as an element in the semiring $\mathbb{Z}^{+}[\omega
_{1}^{\pm },\cdots ,\omega _{m}^{\pm }]$ of Laurent polynomials, the obvious
relations $\omega _{k}^{-1}=$ $\rho \prod\limits_{j\neq k}\omega _{j}$, $%
1\leq k\leq m$, from (D) are sufficient to convert it to an element in $%
\mathbb{Z}^{+}[\omega _{1},\cdots ,\omega _{m},\rho ]$.

\begin{center}
\textbf{4. Preliminaries from representation theory}
\end{center}

Denote by $R^{+}(G)$ the semiring of isomorphism classes $\{V\}$ of $G$%
--Modules $V$ over $\mathbb{C}$. As in \S 3 let $\widetilde{G}$ be the
universal cover of a semi--simple $G$ with $\widetilde{T}\subset \widetilde{G%
}$ the maximal torus that corresponds to $T$ in $G$. The \textsl{character
homomorphism} of $\widetilde{G\text{ }}$ is the semiring map

\begin{center}
$\chi :R^{+}(\widetilde{G})\rightarrow R^{+}(\widetilde{T})=\mathbb{Z}%
^{+}[\omega _{1},\cdots ,\omega _{m},\rho ]$,
\end{center}

\noindent given by restriction $\{V\}\rightarrow \{V\mid \widetilde{T}\}$,
where the identification

\begin{center}
$R^{+}(\widetilde{T})=\mathbb{Z}^{+}[\omega _{1},\cdots ,$ $\omega _{m},\rho
]$
\end{center}

\noindent follows from Lemma 2.

For a $\{V\}\in R^{+}(\widetilde{G})$, the polynomial $\chi (V)\in \mathbb{Z}%
^{+}[\omega _{1},\cdots ,\omega _{m},\rho ]$ is called \textsl{the character
of the }$\widetilde{G}$\textsl{--module }$V$. We will be in particular
interested in the subset $\Omega (G)=$Im$\chi $ of $R^{+}(\widetilde{T})$.
It admits the partition

\begin{center}
$\Omega (G)=\coprod\limits_{1\leq n}$ $\Omega _{n}(G)$
\end{center}

\noindent with $\Omega _{n}(G)=\{g\in \Omega (G)\mid g(1,\cdots ,1)=n\}$.
Alternatively, the $\Omega _{n}(G)$ consists of all characters of $n$%
--dimensional representations of $\widetilde{G}$.

Indeed, elements in $\Omega (G)$ (resp. in $\Omega _{n}(G)$) are abundant.
To explain this we consider the set $\Lambda (G)=\{a_{1}\omega _{1}+\cdots
+a_{m}\omega _{m}\in H^{2}(G/T)\mid a_{m}\in \mathbb{Z}^{+}\}$ of \textsl{%
dominant weights} of $G$. For a $\lambda \in \Lambda (G)$ let $L_{\lambda }$
be the complex line bundle over $G/T$ with first Chern class $%
c_{1}(L_{\lambda })=\lambda $, and let $V_{\lambda }$ be the $\widetilde{G}$%
--module of holomorphic sections of $L_{\lambda }$. Since the set $%
\{V_{\lambda }\mid \lambda \in \Lambda (G)\}$ presents all finite
dimensional complex representations of $\widetilde{G}$ (cf. [LG, p. 69]) and
since the $\chi $ is injective, we have

\bigskip

\textbf{Lemma 4.} \textsl{The set of infinitely many polynomials}

\begin{center}
$\{\chi (V_{\lambda })\in \mathbb{Z}^{+}[\omega _{1},\cdots ,\omega
_{m},\rho ]\mid \lambda \in \Lambda (G)\}$\textsl{\ }
\end{center}

\noindent \textsl{is a basis for the semiabelian group }$\Omega (G)$\textsl{%
\ (over }$\mathbb{Z}^{+}$\textsl{).}

\bigskip

We emphasis at this stage that there have been several methods computing the
character $\chi (V_{\lambda })\in \mathbb{Z}^{+}[\omega _{1},\cdots ,\omega
_{m},\rho ]$ in terms of the $\lambda \in \Lambda (G)$. These are the Weyl
character formula [Hu], Demazure character formula [D] and the
Lakshmibai-Seshadri character formula ([LS], [Li]). Consequently, the $%
\Omega (G)$ admits concrete presentation as a set of polynomials. We present
such an example for further use.

If $G=SU(m)$ is the special unitary group of order $m$, then $(\widetilde{G},%
\widetilde{T})=(G,T)$, $G/T=F(m)$. Consider the semiring map

\begin{center}
$\alpha :R^{+}(T)=\mathbb{Z}^{+}[\omega _{1},\cdots ,\omega _{m-1},\rho
]\rightarrow \mathbb{Z}^{+}[y_{1},\cdots ,y_{m}]$
\end{center}

\noindent by $\omega _{k}\rightarrow y_{1}\cdots y_{k}$, $1\leq k\leq m-1$; $%
\rho \rightarrow y_{2}y_{3}^{2}\cdots y_{m}^{m-1}$, where $\mathbb{Z}%
^{+}[y_{1},\cdots ,y_{m}]$ is the the semiring of polynomials in the $%
y_{1},\cdots ,y_{m}$ with coefficents in $\mathbb{Z}^{+}$ that is subject to
the relation $y_{1}\cdots y_{m}=1$.

For a partition $\mu =(\mu _{1}\geq \mu _{2}\geq \cdots \geq \mu _{m}\geq 0)$
let $s_{\mu }(y_{1},\cdots ,y_{m})\in \mathbb{Z}^{+}[y_{1},\cdots ,y_{m}]$
be \textsl{Schur symmetric function} associated to $\mu $ [M, p.40], and let
$S[y_{1},\cdots ,y_{m}]\subset \mathbb{Z}^{+}[y_{1},\cdots ,y_{m}]$ be the
semiabelian group generated over $\mathbb{Z}^{+}$ by all Schur functions $%
s_{\mu }(y_{1},\cdots ,y_{m})$ associated to the $\mu $ with last part $\mu
_{m}=0$. Then (cf. [FH, Chapter 15], [LG; Chapter 5])

\bigskip

\textbf{Lemma 5.} $\alpha $\textsl{\ maps }$\Omega (SU(m))$\textsl{\
isomorphically onto }$S[y_{1},\cdots ,y_{m}]$.

\bigskip

For the geometries underlying the transformations $\alpha $, we refer to
[DZZ; Example 3] or [DZ$_{2}$; (4.1)].

\begin{center}
\textbf{5. Character factors through homotopy}
\end{center}

Let $Hom^{0}(\widetilde{G},U(n))$ be the set of all homomorphisms $g:%
\widetilde{G}\rightarrow U(n)$ that satisfy $g(\widetilde{T})\subset T^{n}$.
The obvious inclusion $i_{n}:Hom^{0}(\widetilde{G},U(n))\rightarrow R^{+}(%
\widetilde{G})$ is surjective onto the subset $R_{n}^{+}(\widetilde{G})$ of
isomorphism classes of $n$--dimensional complex representations of $%
\widetilde{G}$. On the other hand, since every $g\in $ $Hom^{0}(\widetilde{G}%
,U(n))$ preserves the maximal torus, there is a ready--made map

\begin{center}
$\rho :Hom^{0}(\widetilde{G},U(n))\rightarrow \lbrack G/T,F(n)]$
\end{center}

\noindent defined by the commutivity of the diagram

\begin{center}
$%
\begin{array}{ccc}
\widetilde{G} & \overset{g}{\rightarrow } & U(n) \\
\downarrow &  & \downarrow \\
G/T=\widetilde{G}/\widetilde{T} & \overset{\rho (g)}{\rightarrow } & F(n)%
\end{array}%
$,
\end{center}

\noindent where the vertical maps are the obvious quotients. We may now
organize the relevant correspondences in the diagram below.

\begin{center}
$%
\begin{array}{ccc}
Hom^{0}(\widetilde{G},U(n)) & \overset{i_{n}}{\rightarrow } & R^{+}(%
\widetilde{G}) \\
\rho \downarrow &  & \downarrow \chi \\
\lbrack G/T,F(n)]\overset{r}{\rightarrow } & Hom[H^{2}(F(n)),H^{2}(G/T)] &
\overset{s}{\rightarrow }\mathbb{Z}^{+}[\omega _{1},\cdots ,\omega _{m},\rho
]%
\end{array}%
$.
\end{center}

\bigskip

\textbf{Theorem }(\textsl{Character factors through homotopy})\textbf{. }$%
\chi \circ i_{n}=s\circ r\circ \rho $.

\bigskip

The proof of the Theorem will be postponed until the next section and at
this moment, we show how it leads to a partial solution to our problem. From
the Theorem we obtain

\textbf{Corollary 1. }\textsl{For any polynomial }$g\in \Omega _{n}(G)$%
\textsl{, there exists} \textsl{a map }$f:G/T\rightarrow F(n)$\textsl{\ such
that }$s(f^{\ast })=g$\textsl{. Furthermore, one may read the induced map }$%
f^{\ast }:H^{2}(F(n))\rightarrow H^{2}(G/T)$\textsl{\ from the polynomial }$%
g $\textsl{.}

\bigskip

The first part of Corollary 1 may be rephrased as

\bigskip

\textbf{Corollary 2.} \textsl{For a }$h\in Hom(H^{2}(F(n)),H^{2}(G/T))$%
\textsl{\ with }$s(h)\in \Omega \mathbb{(}G)$\textsl{,} \textsl{there exists
a map }$f:G/T\rightarrow F(n)$\textsl{\ such that }$f^{\ast }=h$\textsl{.}

\bigskip

In algebraic geometry, Schur functions appear as the polynomial
representatives of Schubert classes in complex Grassmanian manifolds (cf.
[L], [St]). In representation theory, Schur functions occur as the
characters of irreducible representations of unitary groups [T], [FH].
Combining Lemma 5 with Corollary 1 we get the following result illustrating
the fashion by which homotopy theory interacts with combinatorics of
symmetric functions.

\bigskip

\textbf{Corollary 3.} \textsl{For any symmetric function \ }$g\in
S[y_{1},\cdots ,y_{m}]$\textsl{\ with }$g(1,\cdots ,1)$ $=n$\textsl{, there
exists a map }$f:F(m)\rightarrow F(n)$\textsl{\ such that }$\alpha \circ
s(f^{\ast })=g$\textsl{.}

\bigskip

\textbf{Remark 2. }According to Corollary 3, every Schur function $s_{\mu
}(y_{1},\cdots ,y_{m})$ with $s_{\mu }(1,\cdots ,$ $1)=n$ \ can arise as the
$\alpha \circ s$--image of a linear map

\begin{center}
$f_{\mu }^{\ast }:H^{2}(F(n))\rightarrow H^{2}(F(m))$.
\end{center}

\noindent We refer to [P] for three definitions of Schur functions in
combinatorics. It would be natural to expect that useful properties of Schur
functions can be derived from the linear maps $f_{\mu }^{\ast }$.

\begin{center}
\textbf{6. Proof of the Theorem}
\end{center}

The binary operation

\begin{center}
$Hom(T,S^{1})\times Hom(T,S^{1})\rightarrow Hom(T,S^{1})$
\end{center}

\noindent by $\alpha \beta =\mu \circ (\alpha \times \beta )\circ \Delta $
furnishes the set $Hom(T,S^{1})$ with the structure of an abelian group,
where $\Delta :T\rightarrow T\times T$ is the diagonal embedding and where $%
\mu :S^{1}\times S^{1}\rightarrow S^{1}$ is the product in $S^{1}$. Consider
the map

\begin{center}
$\iota :Hom(T,S^{1})\rightarrow H^{1}(T)$
\end{center}

\noindent by $\iota (\alpha )=\alpha ^{\ast }[S^{1}]$, where $\alpha ^{\ast
}:H^{1}(S^{1})\rightarrow H^{1}(T)$ is induced by $\alpha $ and where $%
[S^{1}]\in H^{1}(S^{1})=\mathbb{Z}$ is the fundamental class of the oriented
circle $S^{1}$. The following standard fact clarifies the identification $%
H^{1}(T)=Hom(T,S^{1})$ that we have adopted in (B) and in Lemma 2.

\textbf{Lemma 6.} \textsl{The correspondence }$\iota $\textsl{\ is an
isomorphism of abelian groups. In particular, one has}

\begin{center}
$\iota (\omega _{1}^{a_{1}}\cdots \omega _{n}^{a_{n}})$\textsl{\ }$%
=a_{1}\omega _{1}+\cdots +a_{m}\omega _{m}$\textsl{\ (in }$H^{1}(T)$\textsl{%
),}
\end{center}

\noindent \textsl{where we use }$\omega _{i}\in H^{1}(T)$\textsl{\ instead
of }$\iota (\omega _{i})$\textsl{\ (as in Lemma 2 for the sake of
simplicity).}

\bigskip

We are now ready to establish the Theorem. Each representation $g\in $ $%
Hom^{0}(\widetilde{G},$ $U(n))$ induces a bundle map

\begin{center}
$%
\begin{array}{ccccc}
\widetilde{T} & \hookrightarrow & \widetilde{G} & \rightarrow & G/T \\
g^{\prime }\downarrow \quad &  & g\downarrow \quad &  & \rho (g)\downarrow
\quad \\
T^{n} & \hookrightarrow & U(n) & \rightarrow & F(n)%
\end{array}%
,$
\end{center}

\noindent where the $g^{\prime }$ is the restriction of $g$ to $\widetilde{T}
$. From the naturality of the transgression (B) we have the commutative
diagram

\begin{center}
$%
\begin{array}{ccc}
H^{1}(T^{n}) & \overset{\tau _{U(n)}}{\rightarrow } & H^{2}(F(n)) \\
g^{\prime \ast }\downarrow \quad &  & \rho (g)^{\ast }\downarrow \quad \quad
\\
H^{1}(\widetilde{T}) & \underset{\cong }{\overset{\tau _{\widetilde{G}}}{%
\rightarrow }} & H^{2}(G/T)%
\end{array}%
$.
\end{center}

\noindent With the $\tau _{U(n)}$ and $\tau _{\widetilde{G}}$ being
specified respectively in Lemma 1 and 2, we find that if the $\rho (g)^{\ast
}$ is given by

\begin{center}
\noindent $\rho (g)^{\ast }(t_{k})=a_{k,1}\omega _{1}+\cdots +a_{k,m}\omega
_{m}$, $1\leq k\leq n-1$ (cf. (C)),
\end{center}

\noindent in $H^{2}(G/T)$, we have

\begin{enumerate}
\item[(E)] $\qquad g^{\prime \ast }(t_{k}^{\prime })=a_{k,1}\omega
_{1}+\cdots +a_{k,m}\omega _{m}$, $1\leq k\leq n-1$, and

$\qquad g^{\prime \ast }(t_{n}^{\prime })=b_{1}\omega _{1}+\cdots
+b_{m}\omega _{m}$
\end{enumerate}

\noindent in $H^{1}(\widetilde{T})$, where $b_{i}=-(a_{1,i}+\cdots
+a_{n-1,i})$. However, as is standard, the character $\chi (g)$ can be
computed in terms of (E) as

\begin{center}
$\chi (g)=\iota ^{-1}[g^{\prime \ast }(t_{1}^{\prime })]+\cdots +\iota
^{-1}[g^{\prime \ast }(t_{n}^{\prime })]$
\end{center}

\noindent This agrees with $s(\rho (g)^{\ast })$ by Lemma 6.

\begin{center}
\textbf{7. Endnotes}
\end{center}

This paper is by no means a final exposition in the topic, but is part of a
larger effort to understand maps between flag manifolds by means of
cohomology.

\bigskip\

\textbf{7.1. }During the past two decades many works were devoted to the
study of self--maps of flag manifolds, cf. [GH], [H$_{1}$], [H$_{2}$], [HH],
[P$_{1}$], [P$_{2}$], [DZ$_{1}$], [Du]. However, the corresponding
investigation into maps between different flag manifolds has not yet
received as much attention as it deserves. The aim of present work is to
demonstrate concrete examples (from representation theory) indicating the
richness of the latter topic.

\bigskip

\textbf{7.2.} Our method can be extended to the study of homotopy classes of
maps $G/T\rightarrow G^{\prime }/T^{\prime }$ between two arbitray flag
manifolds. For instance in the Problem in \S 2, one may replace the flag
manifold $F(n)$ by

\begin{center}
$D(n)=SO(n)/T$ (resp. $T(n)=Sp(n)/T$ ),
\end{center}

\noindent where the $SO(n)$ (resp. $Sp(n)$) is the special orthogonal group
(resp. symplectic group) of order $n$ and where $T\subset SO(n)$ (resp. $%
T\subset Sp(n)$) is a maximal torus. In these cases results analogous to the
Theorem hold, only the character homomorphisms for the real and quaternionic
representations will be involved and the correspondence $s$ in \S 3 would be
modified accordingly.

\bigskip

\textbf{7.3.} Let $Map(G/T,F(n))$ be the space of all continuous maps $%
G/T\rightarrow F(n)$. For a $f\in Map(G/T,F(n))$, let $M_{f}\subset
Map(G/T,F(n))$ be the path--connected component containing $f$. In [Me] and
[S], W. Meier and S. Smith initiated the project to determine the rational
homotopy type of the space $M_{f}$.

\begin{quote}
\textbf{Conjecture. }The rational homotopy type of $M_{f}$ is determined by
the polynomial $s\circ r(f)\in \mathbb{Z}^{+}[\omega _{1},\cdots ,\omega
_{m},\rho ]$.
\end{quote}

\bigskip

More precisely, if $g\in Hom^{0}(U(m),U(n))$ is such that

\begin{center}
$\alpha \circ \chi (g)=l(t_{1}+\cdots +t_{m})+(n-lm)$ (cf. Corollary 3),
\end{center}

\noindent then, modulo the products of odd dimensional spheres, the rational
homotopy type of the space $M_{\rho (g)}$ has been determined by S. Smith
[S]. We observe that the sum $t_{1}+\cdots +t_{m}$ is the Schur function
associated to the partition $\mu =(1,0,\cdots ,0)$ ([M, P.42]). We note also
from Lemma 5 that corresponding to any Schur function $s_{\mu }(y_{1},\cdots
,y_{m})$ with $s_{\mu }(1,\cdots ,1)=n$, there is a representation $g_{\mu
}\in Hom^{0}(U(m),U(n))$ with $\alpha \circ \chi (g_{\mu })=s_{\mu
}(y_{1},\cdots ,y_{m})$.

\begin{quote}
\textbf{Problem.} For a partition $\mu =(\mu _{1},\mu _{2},\cdots ,\mu _{m})$%
, determine the rational homotopy type of the space $M_{\rho (g_{\mu })}$.
\end{quote}

\bigskip

\textbf{Acknowledgements.} This paper is based on a talk given at the
International conference on\textsl{\ homotopy theory and related topics}
held in Seoul, February 1-4, 2005. The author is grateful to S. Smith for
suggestive discussion during the conference.

\begin{center}
\textbf{References}
\end{center}

\begin{enumerate}
\item[{[B]}] A. Borel, Sur la cohomologie des espaces fibres principaux et
des homogenes de groupes de Lie compacts, Ann. Math., 57(1953), 115-207.

\item[{[BH]}] A. Borel and F. Hirzebruch, Characteristic classes and
homogeneous space, I., Amer. J. Math., 80 (1958), 458-538.

\item[{[BS]}] R. Bott and H. Samelson, The cohomology ring of $G/T$, Nat.
Acad. Sci. 41 (7) (1955), 490-492.

\item[{[BGG]}] N. Bernstein, I.M. Gelfand and S. I. Gelfand, Schubert cells
and cohomology of the spaces $G/P$, Russian Math. Surv. 28(1973), 1-26.

\item[{[D]}] M. Demazure, D\'{e}singularization des vari\'{e}t\'{e}s de
Schubert g\'{e}n\'{e}ralis\'{e}es, Ann. Sci. \'{E}cole. Norm. Sup. (4)
7(1974), 53-88.

\item[{[Du]}] H. Duan, Self-maps of the Grassmannian of complex structures.
Compositio Math. 132 (2002), no. 2, 159--175.

\item[{[DZ$_{1}$]}] H. Duan and Xu-an Zhao, The classification of cohomology
endomorphisms of certain flag manifolds. Pacific J. Math. 192 (2000), no. 1,
93--102.

\item[{[DZ$_{2}$]}] H. Duan, Xu--an Zhao, The height function on the
2-diemensional cohomology of a flag manifold, J. Lie Theory, 15(1)(2005),
219-226.

\item[{[DZZ]}] H. Duan, Xu--an Zhao and Xuezhi Zhao, The Cartan matrix and
enumerative calculus. J. Symbolic Comput. 38 (2004), no. 3, 1119--1144.

\item[{[FH]}] W. Fulton and J. Harris, Representation theory. A first course.
Graduate Texts in Mathematics, 129. Springer-Verlag, New York, 1991.

\item[{[GH]}] H. Glover,and W. Homer, Self-maps of flag manifolds. Trans.
Amer. Math. Soc. 267 (1981), no. 2, 423--434.

\item[{[H$_{1}$]}] M. Hoffman, On fixed point free maps of the complex flag
manifold, Indiana Univ. Math. J. 33 (1984), no. 2, 249--255.

\item[{[H$_{2}$]}] M. Hoffman, Endomorphisms of the cohomology of complex
Grassmannians. Trans. Amer. Math. Soc. 281 (1984), no. 2, 745--760.

\item[{[HH]}] M. Hoffman and W. Homer, On cohomology automorphisms of complex
flag manifolds, Proc. AMS. 91(4)(1984), 643-648.

\item[{[Hu]}] J. E. Humphreys, Introduction to Lie algebras and
representation theory, Graduated Texts in Math. 9, Springer-Verlag New York,
1972.

\item[{[K]}] S. Kumar, Kac-Moody groups, their flag varieties and
representation theory. Progress in Mathematics, 204. Birkh\"{a}user Boston,
Inc., Boston, MA, 2002. xvi+606 pp.

\item[{[L]}] L. Lesieur, Les problemes d'intersections sur une variete de
Grassmann, C. R. Acad. Sci. Paris, 225 (1947), 916-917.

\item[{[Li]}] P. Littelmann, A Littelwood--Richardson rule for symmetrizable
Kac--Moody algebra, Invent. Math. 116(1994), 329-346.

\item[{[LG]}] V. Lakshmibai and N. Gonciulea, Flag varieties,
Hermann-Acutalities Mathematiques, Hermann Editeurs Des Sciences et Des
Arts, 2001.

\item[{[LS]}] V. Lakshmibai and C. S. Seshadri, Standard monomial theory and
Schubert varieties--a survey, Proceedings of the conference on ''Algebraic
groups and Applications'', 1991, 279-322, Manoj Prakashan.

\item[{[M]}] I. G. Macdonald, Symmetric functions and Hall polynomials,
Oxford Mathematical Monographs, Oxford University Press, Oxford, second ed.,
1995.

\item[{[Me]}] W. Meier, Rational universal fibrations and flag manifolds,
Math. Ann. 258(1982), 329--340.

\item[{[P]}] P. Pragacz, Algebro-geometric applications of Schur S-- and
Q--polynomials, Topics in invariant Theory (M.-P. Malliavin, ed.), Lecture
Notes in Math., Vol. 1478, Spring-Verlag, Berlin and New York, 1991, 130-191.

\item[{[P$_{1}$]}] S. Papadima, Rigidity properties of compact Lie group
module maximal tori, Math. Ann. 275(1986) 637-652.

\item[{[P$_{2}]$}] S. Papadima, Rational homotopy equivalences of Lie type.
Math. Proc. Cambridge Philos. Soc. 104 (1988), no. 1, 65--80.

\item[{[S]}] S. Smith, Rational classification of simple function space
components for flag manifolds. Canad. J. Math. 49 (1997), no. 4, 855--864.

\item[{[St]}] R. Stanley, Some combinatorial aspects of the Schubert
calculus, Combinatoire et repr\'{e}sentation du groupe sym\'{e}trique,
Strasbourg (1976), 217-251.

\item[{[T]}] H. Tamvakis, The connection between representation theory and
Schubert calculus, \textsl{preprint available on} arXiv: math.AG/0306414.
\end{enumerate}

\end{document}